\documentclass[EJP,preprint]{ejpecp}
\usepackage[T1]{fontenc}
\usepackage{enumitem}
\hypersetup{unicode=true,
  pdftitle={Ultra-hard-edge universality under deterministic preconditioning},
  pdfauthor={Christian E. H\"aggblom}}

\SHORTTITLE{Anisotropic ultra-hard-edge universality}
\TITLE{Ultra-hard-edge universality under deterministic preconditioning}
\AUTHORS{Christian E. H\"aggblom}
\AMSSUBJ{60B20}

\newcommand{\R}{\mathbb R}
\newcommand{\E}{\mathbb E}
\newcommand{\Prob}{\mathbb P}
\newcommand{\1}{\mathbf 1}
\newcommand{\Tr}{\operatorname{Tr}}
\newcommand{\rank}{\operatorname{rank}}
\newcommand{\diag}{\operatorname{diag}}
\newcommand{\op}{\mathrm{op}}
\newcommand{\HS}{\mathrm{HS}}
\newcommand{\Haar}{\mathrm{Haar}}
\newcommand{\Gn}{\mathcal G_n}
\newcommand{\good}{\mathcal Q_n}
\newcommand{\reg}{\mathcal R_n}
\newcommand{\smin}{s_{\min}}
\newcommand{\norm}[1]{\left\lVert #1\right\rVert}

\newcommand{\ip}[2]{\langle #1,#2\rangle}
\newcommand{\law}{\mathcal L}

\ABSTRACT{
Let $A_n$ be an $n\times n$ real matrix with iid entries distributed as a fixed centered, variance-one, symmetric subgaussian random variable. For deterministic invertible $M_n$, set $h_n=\|M_n^{-1}\|_{\mathrm{HS}}$ and $B_n=M_n^{-1}/h_n$. We prove
\[
\mathbb{P}(s_{\min}(M_nA_n)\le\varepsilon/h_n)
=(1+o(1))\varepsilon\Theta_n(B_n)
\]
as $n\to\infty$, uniformly over $M_n$ and $e^{-cn}\le\varepsilon\le n^{-C}$. Here $\Theta_n(B)=\mathbb{E}\sqrt{nU^\top BB^\top U}$, where $U$ is uniform on $S^{n-1}$. The constants $c>0$ and $C>1$ depend only on the entry law. No condition-number bound or effective-rank assumption is required.

The coefficient tends to one for diffuse inverse profiles. For convergent concentrated or spiked profiles, its limit is an explicit Gaussian expectation. The proof combines the submicroscopic comparison of Sah, Sahasrabudhe, and Sawhney with signed-permutation averaging and a weighted comparison for random hyperplane normals. Applications include inverse amplification, square sample covariance matrices with arbitrary deterministic positive-definite population matrices, and independent random environments.
}

\begin{document}

\section{Introduction and main results}

\subsection{The anisotropic lower tail}
For a square iid random matrix, multiplication by a deterministic matrix stretches different directions by different amounts. We ask how this changes the probability that the resulting product is exceptionally close to singular. Our result separates the effect of this anisotropy into an overall scale, determined by the Hilbert--Schmidt norm of the inverse multiplier, and a leading coefficient given by a spherical average. The coefficient records whether the inverse multiplier acts comparably in many directions or concentrates its action in a few, while remaining universal across the entry distributions considered here.

Throughout the paper, $\xi$ is a fixed real random variable satisfying
\begin{equation}\label{eq:entry}
\E\xi=0,\qquad \E\xi^2=1,\qquad \xi\overset{d}=-\xi,
\qquad \norm{\xi}_{\psi_2}<\infty.
\end{equation}
Let $A_n$ be an $n\times n$ matrix whose entries are independent copies of $\xi$. For a deterministic invertible matrix $M$, set
\begin{equation}\label{eq:normalization}
h(M)=\norm{M^{-1}}_{\HS},\qquad
B=\frac{M^{-1}}{h(M)},\qquad H=BB^\top,\qquad \Tr H=1.
\end{equation}
Here $\norm T_{\HS}^2=\Tr(TT^\top)$ is the squared Hilbert--Schmidt norm.
For $H\succeq0$ with $\Tr H=1$ and $B\in\R^{n\times m}$ with $\norm B_{\HS}=1$, define
\begin{equation}\label{eq:theta}
\vartheta_n(H):=\E_{U\sim\Haar(S^{n-1})}\sqrt{nU^\top HU},
\qquad \Theta_n(B):=\vartheta_n(BB^\top).
\end{equation}
Thus $\Theta_n(B)$ depends only on the eigenvalues of $BB^\top$. For any nonzero deterministic $K\in\R^{n\times m}$,
\[
\Theta_n\!\left(\frac{K}{\norm K_{\HS}}\right)
=\frac{\sqrt n}{\norm K_{\HS}}\,\E_U\norm{K^\top U}_2.
\]

For $M=I$, one has $h(I)=\sqrt n$, $H=I/n$, and $\Theta_n(I/\sqrt n)=1$. The ordinary hard-edge scale is then $n^{-1/2}$. We use the term \emph{ultra-hard edge} for the smaller thresholds $\varepsilon/h(M)$ considered below, where $\varepsilon$ lies between an exponentially small lower cutoff and an inverse-polynomial upper cutoff. The theorem determines the leading coefficient throughout this shrinking range, uniformly over the preconditioner.

\begin{theorem}[Deterministic preconditioning]\label{thm:main}
Assume \eqref{eq:entry}. There are constants $0<c_0<c_1$, $C_0>1$, $C<\infty$, and a deterministic sequence $r_n\to0$, depending only on the entry law, such that for every deterministic invertible $n\times n$ matrix $M$ and every
\begin{equation}\label{eq:window}
e^{-c_0n}\le\varepsilon\le n^{-C_0},
\end{equation}
with $B$ defined in \eqref{eq:normalization},
\begin{equation}\label{eq:main-quant}
\left|\Prob\!\left(\smin(MA_n)\le
\frac{\varepsilon}{\norm{M^{-1}}_{\HS}}\right)
-\varepsilon\Theta_n(B)\right|
\le C\varepsilon\left[r_n+
\bigl(\varepsilon n\norm B_{\op}^2\bigr)^{1/3}\right]
+Ce^{-c_1n}.
\end{equation}
All bounds are uniform over $M$. In particular, throughout \eqref{eq:window},
\begin{equation}\label{eq:main-relative}
\Prob\!\left(\smin(MA_n)\le
\frac{\varepsilon}{\norm{M^{-1}}_{\HS}}\right)
=\bigl(1+o(1)\bigr)\varepsilon\Theta_n(B).
\end{equation}
\end{theorem}

The coefficient is uniformly nondegenerate:
\begin{equation}\label{eq:theta-bounds-intro}
\sqrt{\frac2\pi}\le\Theta_n(B)\le1.
\end{equation}
We prove this and more precise profile statements in Section~\ref{sec:profiles}. Since $\norm B_{\op}\le1$, the explicit residual in \eqref{eq:main-quant} is uniformly negligible in \eqref{eq:window}. The inequality $c_0<c_1$ ensures that the additive exponential remainder is negligible relative to every threshold in that window. The theorem makes no assertion about the optimal exponential endpoint.

\begin{corollary}[Gaussian-reference universality]\label{cor:gaussian}
Let $G_n$ have iid standard real Gaussian entries. Under the assumptions and in a common window of the form \eqref{eq:window},
\begin{equation}\label{eq:gaussian-reference}
\frac{\Prob(\smin(MA_n)\le\varepsilon/\norm{M^{-1}}_{\HS})}
{\Prob(\smin(MG_n)\le\varepsilon/\norm{M^{-1}}_{\HS})}
\longrightarrow1
\end{equation}
uniformly over deterministic invertible $M$ and admissible $\varepsilon$.
\end{corollary}

\begin{proof}
Apply Theorem~\ref{thm:main} to both entry laws, shrink the window to their intersection, and use \eqref{eq:theta-bounds-intro}.
\end{proof}

\subsection{Diffuse and concentrated inverse profiles}
Define
\begin{equation}\label{eq:diffuse}
\Delta_4(M)^2=\Tr\bigl((BB^\top)^2\bigr).
\end{equation}
The quantity $\Delta_4(M)$ tends to zero exactly when the inverse singular mass is diffuse. Theorem~\ref{thm:main} includes this regime, but does not require it.

\begin{corollary}[Diffuse coefficient]\label{cor:diffuse}
If $\Delta_4(M_n)\to0$, then uniformly in \eqref{eq:window},
\[
\Prob\!\left(\smin(M_nA_n)\le
\frac{\varepsilon}{\norm{M_n^{-1}}_{\HS}}\right)
=(1+o(1))\varepsilon.
\]
More precisely, \eqref{eq:main-quant} remains valid with $\varepsilon\Theta_n(B)$ replaced by $\varepsilon$, provided $\varepsilon\Delta_4(M_n)^2$ is added to its right-hand side.
\end{corollary}

At the other extreme, when the normalized inverse singular mass approaches rank one, the coefficient tends to $\sqrt{2/\pi}$. More generally, if finitely many weights $\lambda_1,\ldots,\lambda_r$ survive and the residual mass $d=1-\sum_{j=1}^r\lambda_j$ is diffuse, the coefficient tends to
\begin{equation}\label{eq:spike-intro}
\E\sqrt{d+\sum_{j=1}^r\lambda_jg_j^2},
\qquad g_j\ \text{independent }N(0,1).
\end{equation}
Thus deterministic anisotropy survives through a Gaussian coefficient even though the entry distribution disappears from the leading term.

\subsection{Related work}
Edelman computed the limiting distribution of the least singular value in the Gaussian case \cite{Edelman}. Tao and Vu established universality of this distribution for iid matrices under a sufficiently high finite-moment assumption, with a polynomially small additive error uniform over thresholds \cite{TV}. For correlated complex Gaussian Wishart matrices with a fixed difference between the matrix dimensions and population eigenvalues bounded above and away from zero, Hachem, Hardy, and Najim obtained Bessel-kernel limits at the hard edge \cite{HHN}. These results describe distributions on the ordinary hard-edge scale. Relative asymptotics for exponentially small lower-tail probabilities require error estimates that remain negligible compared with those probabilities. A polynomially small additive error alone is insufficient.

Rudelson and Vershynin established linear lower-tail bounds for square iid subgaussian matrices and the corresponding optimal rectangular estimates \cite{RVSquare,RVRect}. Sah, Sahasrabudhe, and Sawhney (SSS) proved the sharp unpreconditioned coefficient at exponentially small relative thresholds \cite{SSS}. Letwin and Polavarapu identified $\norm{M^{-1}}_{\HS}^{-1}$ as the general preconditioned scale under substantially weaker entry assumptions \cite{LP}. Their theorem supplies an upper bound with an unspecified constant. The present result identifies the coefficient under \eqref{eq:entry}, uniformly over deterministic invertible preconditioners, without a condition-number bound or an effective-rank assumption. The covariance consequence in Corollary~\ref{cor:covariance} concerns real iid latent entries satisfying \eqref{eq:entry} and allows arbitrary positive-definite population matrices, with no uniform spectral bounds, at the ultra-small thresholds specified there.

The inverse-tail mechanism used here originated in the Gaussian phase-retrieval argument at the critical injectivity threshold \cite[Section~9.1, Lemma~7]{HaggblomPR}. There, conditioning on the common rows of overlapping square minors leads to a weighted inverse-tail problem for an independent square Gaussian block. Its leading coefficient involves the spherical average $\E_U\norm{K^\top U}_2$, where $U$ is Haar-uniform on $S^{n-1}$. For a nonzero $n\times m$ input matrix $K$, this average satisfies
\[
\E_U\norm{K^\top U}_2
=\frac{\norm K_{\HS}}{\sqrt n}\,
\Theta_n\!\left(\frac{K}{\norm K_{\HS}}\right).
\]
Thus the same coefficient functional already occurs in that Gaussian argument. The earlier inverse-tail estimate uses the Haar distribution of Gaussian singular vectors and their independence from the singular values, together with moment assumptions on the random input that force its normalized profile to be diffuse. The present paper establishes universality of this coefficient for symmetric iid subgaussian entries, uniformly over arbitrary deterministic input profiles, including concentrated ones. The signed-permutation arguments below replace the exact rotational invariance used in the Gaussian analysis.

\subsection{Proof strategy}
The proof separates the rare-event estimate from the identification of its coefficient. A one-column inverse decomposition produces the mark
\[
a_B(C)=\sqrt n\,\norm{B^\top v(C)}_2,
\]
where $C$ is obtained by deleting one column of $A_n$ and $v(C)$ is a unit normal to its column space. This mark measures the inverse multiplier in the random normal direction. It need not concentrate when the inverse profile has a few large eigenvalues.

The first ingredient is the marked transfer theorem, Theorem~\ref{thm:transfer}. It extends the SSS unmarked regularity and conditional replacement estimates to marks whose second moments are bounded on every signed-permutation orbit. For $a_B$, averaging $a_B(PC)^2$ over signed permutation matrices $P$ gives exactly one. This identity permits truncation at $a_B\le L$ with rare-event cost $O(\varepsilon/L)$. On the retained set, the unmarked estimates apply at the enlarged threshold $\varepsilon L$. The resulting coefficient is $\E[a_B(C)\Lambda_n(C)]$, where $\Lambda_n$ is a signed-permutation-invariant spectral intensity.

The second ingredient is the weighted coefficient comparison, Theorem~\ref{thm:weighted}. It identifies this expectation uniformly in $B$, using an $L^p$ bound for $\Lambda_n$ with $1<p<2$ and weak convergence of two coordinate magnitudes of $v(C)$. No independence between the normal and the intensity is required. The comparison follows from a transport inequality for signed-permutation averages. The coordinate limit is established in Lemma~\ref{lem:coordinate} and Appendix~\ref{app:coordinate} using the projected-column method of Nguyen and Vu \cite{NV}, with an explicit finite-population variance estimate for the sampling normalization. Combining the two ingredients gives the universal coefficient $\Theta_n(B)$.

\paragraph{Organization.}
Section~\ref{sec:geometry} gives the inverse decomposition. Sections~\ref{sec:orbits}--\ref{sec:transfer} develop the weighted comparison, the spectral intensity, and the marked transfer theorem. Section~\ref{sec:main-proof} combines these ingredients to prove Theorem~\ref{thm:main}. Sections~\ref{sec:profiles} and~\ref{sec:applications} describe the coefficient and its applications. Sections~\ref{sec:scope} and~\ref{sec:directions} discuss the scope of the method and formulate extensions, with exact Gaussian identities. Appendix~\ref{app:coordinate} proves the coordinate limit.

\subsection{Conventions}
We write $\norm T_{\op}$ for the Euclidean operator norm and use the subgaussian norm
\[
\norm\xi_{\psi_2}:=\inf\{s>0:\E\exp(\xi^2/s^2)\le2\}.
\]
Constants may depend on the fixed entry law. The notation $o(1)$ always denotes a deterministic quantity tending to zero. When matrices or thresholds vary, uniformity is stated explicitly. Singular values of a deleted $n\times(n-1)$ matrix are ordered increasingly as $s_1(C)\le\cdots\le s_{n-1}(C)$. Square least singular values are denoted by $\smin$. All assertions concerning a normal vector are insensitive to its global sign. Rank-deficient deleted matrices are excluded from geometric events and assigned intensity zero. Their probability is $O(e^{-cn})$. In inverse-amplification statements, we count a singular $A_n$ as a failure and assign its inverse norm the value $+\infty$.

The logarithmic cutoffs and regularity events below are defined for $n\ge n_0$, with $n_0$ sufficiently large. Uniform moment bounds for the associated intensities are understood over this range. The value of $n_0$ may be increased in the asymptotic arguments.

\section{The inverse geometry}\label{sec:geometry}

Write $A=[C\ X]$, where $C\in\R^{n\times(n-1)}$ and $X\in\R^n$. On $\{\rank C=n-1\}$ choose $Q$ with orthonormal columns spanning $\operatorname{col}C$, and choose a measurable unit vector $v\in\ker C^\top$. Put
\[
R=Q^\top C,\qquad g=Q^\top X,\qquad x=v^\top X.
\]
Then $(Q,v)^\top A=\left(\begin{smallmatrix}R&g\\0&x\end{smallmatrix}\right)$.

\begin{proposition}[One-column decomposition]\label{prop:inverse}
For any deterministic $B\in\R^{n\times m}$, on $\{\rank C=n-1,x\ne0\}$,
\begin{equation}\label{eq:inverse-decomp}
A^{-1}B=E_C+x^{-1}w_C(X)v^\top B,
\quad
E_C=\begin{pmatrix}R^{-1}Q^\top B\\0\end{pmatrix},
\quad
w_C(X)=\begin{pmatrix}-R^{-1}Q^\top X\\1\end{pmatrix}.
\end{equation}
If
\begin{equation}\label{eq:full-correction}
\Xi_C(X)=\sqrt{1+\norm{R^{-1}Q^\top X}_2^2},\qquad
D_C=\norm{E_C}_{\op},
\end{equation}
then
\begin{equation}\label{eq:inverse-sandwich}
\left|\norm{A^{-1}B}_{\op}
-\frac{\Xi_C(X)\norm{B^\top v}_2}{|v^\top X|}\right|
\le D_C\le\frac{\norm B_{\op}}{s_1(C)}.
\end{equation}
\end{proposition}

\begin{proof}
Block inversion proves \eqref{eq:inverse-decomp}. The second term has rank at most one, with norm
\[
|x|^{-1}\norm{w_C(X)}_2\norm{B^\top v}_2.
\]
Since $C=QR$, the singular values of $R$ are the nonzero singular values of $C$. The triangle inequality gives the remaining assertions. Also $R^{-1}Q^\top=C^\dagger$, so $\Xi_C(X)=\sqrt{1+\norm{C^\dagger X}_2^2}$ is independent of the choice of $Q$.
\end{proof}

For subsequent unconditional probability estimates, set $D_C=+\infty$ on $\{\rank C<n-1\}$.

For a nonnegative $C$-measurable mark $a$, let
\begin{equation}\label{eq:T-def}
T_{ta}=\{\rank C=n-1,\ |v^\top X|\le ta(C)n^{-1/2}\Xi_C(X)\},
\qquad T_s=T_{s\cdot1}.
\end{equation}
Also put $S_s=\{\smin(A)\le s/\sqrt n\}$.

\begin{lemma}[Exact event inclusion]\label{lem:inclusion}
For every $s\ge0$, $T_s\subset S_s$. If $A$ is iid as in \eqref{eq:entry}, then
\begin{equation}\label{eq:unmarked-linear}
\Prob(T_s)\le Cs+Ce^{-cn},\qquad s\ge0.
\end{equation}
\end{lemma}

\begin{proof}
The vector $w_C(X)$ from \eqref{eq:inverse-decomp} satisfies $Aw_C(X)=(v^\top X)v$ and $\norm{w_C(X)}_2=\Xi_C(X)$. Hence $\smin(A)\le |v^\top X|/\Xi_C(X)$, including when $v^\top X=0$. The probability bound follows from the square least-singular-value estimate \cite{RVSquare}.
\end{proof}

For $\norm B_{\HS}=1$ define $a_B(C)=\sqrt n\norm{B^\top v(C)}_2$ on full rank, and zero otherwise. Proposition~\ref{prop:inverse} gives, for $0<\delta<1$, outside singularity and $\{D_C>\delta/\varepsilon\}$,
\begin{equation}\label{eq:marked-sandwich}
T_{\varepsilon a_B/(1+\delta)}
\subset\{\norm{A^{-1}B}_{\op}\ge\varepsilon^{-1}\}
\subset T_{\varepsilon a_B/(1-\delta)}.
\end{equation}
The rectangular estimate \cite{RVRect} gives
\begin{equation}\label{eq:residual}
\Prob(D_C>\delta/\varepsilon)
\le C\frac{\varepsilon^2 n\norm B_{\op}^2}{\delta^2}+Ce^{-cn}.
\end{equation}

\section{Signed-permutation averaging}\label{sec:orbits}

Let $\Gn$ denote the finite group of signed permutation matrices of order $n$. All averages over $P\in\Gn$ use its uniform probability measure. For fixed $z\in\R^n$,
\begin{equation}\label{eq:orbit-cov}
\E_P Pzz^\top P^\top=\frac{\norm z_2^2}{n}I.
\end{equation}
Under \eqref{eq:entry}, $(PC,PX)$ has the same law as $(C,X)$, while the normal line transforms as $v(PC)v(PC)^\top=Pv(C)v(C)^\top P^\top$.

\subsection{Large marks inside an invariant rare event}
The following elementary lemma is the marked-regularity input.

\begin{lemma}[Orbit truncation]\label{lem:orbit-tail}
Let a finite group act on a random object $Z$ with invariant law. Let $(F_s)_{s\ge0}$ be increasing invariant events satisfying
\[
\Prob(F_s)\le Ks+b,\qquad s\ge0.
\]
Suppose $a(Z)\ge0$ satisfies $\E_Pa(PZ)^2\le K_a$ for every $Z$. Then, for $t>0$ and $L\ge1$,
\begin{equation}\label{eq:orbit-tail}
\Prob(F_{ta(Z)}\cap\{a(Z)>L\})
\le CKK_a\frac tL+CK_a\frac b{L^2}.
\end{equation}
If $W(Z)\ge0$ is invariant and integrable, then
\begin{equation}\label{eq:weighted-tail}
\E[aW\1_{\{a>L\}}]\le\frac{K_a}{L}\E W.
\end{equation}
\end{lemma}

\begin{proof}
Split $\{a>L\}$ into the intervals $(2^jL,2^{j+1}L]$, $j\ge0$. Monotonicity, invariance of $F_s$, and Markov's inequality on each orbit give
\begin{align*}
\Prob(F_{ta}\cap\{2^jL<a\le2^{j+1}L\})
&\le\E\left[\1_{F_{t2^{j+1}L}}\E_P\1_{\{a(PZ)>2^jL\}}\right]\\
&\le\frac{K_a}{2^{2j}L^2}\Prob(F_{t2^{j+1}L}).
\end{align*}
Summation gives \eqref{eq:orbit-tail}. For the second assertion, average over the group and use $a\1_{\{a>L\}}\le a^2/L$.
\end{proof}

For the norm mark, \eqref{eq:orbit-cov} gives the exact identity
\begin{equation}\label{eq:mark-second}
\E_Pa_B(PC)^2=1,\qquad \rank C=n-1,
\end{equation}
and consequently $\E a_B\le1$. The full events $T_s$ are invariant under simultaneous signed permutations of $C$ and $X$. Lemmas~\ref{lem:inclusion} and \ref{lem:orbit-tail} therefore imply
\begin{equation}\label{eq:large-mark-cost}
\Prob(T_{ta_B}\cap\{a_B>L\})
\le C\frac tL+Ce^{-cn}.
\end{equation}
This is a bound on the marked rare event itself. No exceptional set of merely polynomial probability is discarded.

\subsection{Transport of orbit-averaged marks}
For $v\in S^{n-1}$ and $\norm B_{\HS}=1$, put
\[
F_B(v)=\sqrt n\norm{B^\top v}_2,\qquad
\overline F_B(v)=\E_PF_B(Pv),\qquad
\mu_v=\frac1n\sum_{i=1}^n\delta_{\sqrt n|v_i|}.
\]
Write $W_2$ for the quadratic Wasserstein distance on probability measures on $[0,\infty)$.

\begin{lemma}[Orbit transport]\label{lem:transport}
For all $u,v\in S^{n-1}$ and all rectangular $B$ with $\norm B_{\HS}=1$,
\begin{equation}\label{eq:transport}
|\overline F_B(v)-\overline F_B(u)|\le W_2(\mu_v,\mu_u).
\end{equation}
Moreover $0\le\overline F_B(v)\le1$.
\end{lemma}

\begin{proof}
Let $v^\uparrow,u^\uparrow$ be the increasingly sorted coordinate magnitudes. They represent the same signed-permutation orbits as $v,u$. Couple their orbit averages using the same $P$, and put $d=\sqrt n(v^\uparrow-u^\uparrow)$. The reverse triangle inequality, Cauchy--Schwarz, and \eqref{eq:orbit-cov} give
\[
|\overline F_B(v)-\overline F_B(u)|
\le\E_P\norm{B^\top Pd}_2
\le\left(\frac{\norm d_2^2}{n}\Tr BB^\top\right)^{1/2}.
\]
Increasing matching is the optimal quadratic coupling of the two equally weighted one-dimensional empirical measures, so the last expression is $W_2(\mu_v,\mu_u)$. Finally $\E_PF_B(Pv)^2=1$, and Jensen's inequality gives the upper bound.
\end{proof}

\begin{theorem}[Weighted coefficient comparison]\label{thm:weighted}
Let $(v_n,W_n)$ be a random unit line and a real random variable whose joint law is invariant under $(v_n,W_n)\mapsto(Pv_n,W_n)$, $P\in\Gn$. Suppose
\begin{equation}\label{eq:pair-normal}
(\sqrt n|v_{n,1}|,\sqrt n|v_{n,2}|)
\ \Longrightarrow\ (|g_1|,|g_2|),
\end{equation}
where $g_1,g_2$ are independent standard Gaussians, and suppose $\sup_n\norm{W_n}_{L^p}<\infty$ for some $p>1$. Then
\begin{equation}\label{eq:weighted-comparison}
\sup_{m\ge1}\ \sup_{B\in\R^{n\times m},\,\norm B_{\HS}=1}
\left|\E[F_B(v_n)W_n]-\Theta_n(B)\E W_n\right|\longrightarrow0.
\end{equation}
No independence between $v_n$ and $W_n$ is assumed.
\end{theorem}

\begin{proof}
Let $\nu=\law(|g_1|)$. For bounded continuous $f$, exchangeability and \eqref{eq:pair-normal} show that
\[
\E\int f\,d\mu_{v_n}\to\int f\,d\nu,
\qquad
\operatorname{Var}\!\left(\int f\,d\mu_{v_n}\right)\to0.
\]
Applying this to a countable convergence-determining class gives $\mu_{v_n}\Rightarrow\nu$ in probability. Both measures have second moment exactly one. The characterization of $W_2$ convergence by weak convergence and convergence of second moments therefore yields
\begin{equation}\label{eq:emp-w2}
W_2(\mu_{v_n},\nu)\to0\quad\text{in probability}.
\end{equation}
This implication can also be applied along almost surely convergent subsubsequences. Since the measures are supported on $[0,\infty)$ and have second moment one, their squared $W_2$ distance is at most two. Thus \eqref{eq:emp-w2} holds in every fixed $L^q$, $q<\infty$.

For an independent Haar vector $U_n$, the same assertion follows by writing $U_n=g/\norm g_2$ and applying the law of large numbers. By Lemma~\ref{lem:transport},
\begin{equation}\label{eq:uniform-orbit-error}
\sup_{B}\left|\overline F_B(v_n)-\Theta_n(B)\right|
\le W_2(\mu_{v_n},\nu)+\E W_2(\mu_{U_n},\nu).
\end{equation}
Here $\E\overline F_B(U_n)=\Theta_n(B)$. Joint invariance gives
\[
\E[F_B(v_n)W_n]=\E[\overline F_B(v_n)W_n].
\]
H\"older's inequality with $q=p/(p-1)$ and \eqref{eq:uniform-orbit-error} finish the proof. All bounds are independent of $B$ and of its number of columns.
\end{proof}

\begin{lemma}[Two coordinate magnitudes]\label{lem:coordinate}
Let $C_n$ be an $n\times(n-1)$ matrix with iid centered, variance-one, subgaussian entries. On full rank, let $v_n$ be a unit vector spanning $\ker C_n^\top$, and on rank deficiency extend its line by an independent Haar line. Then
\[
(\sqrt n\,|v_{n,1}|,\sqrt n\,|v_{n,2}|)
\ \Longrightarrow\ (|g_1|,|g_2|),
\]
where $g_1,g_2$ are independent standard Gaussians.
\end{lemma}

\begin{proof}
Appendix~\ref{app:coordinate} gives a proof using coordinate delocalization and the projected-column comparison of Nguyen--Vu \cite{NV}, with an explicit finite-population variance estimate for the sampling normalization. Only weak convergence of two magnitudes is asserted. The exponentially small rank-deficiency probability does not affect this limit.
\end{proof}

Under \eqref{eq:entry}, the Haar-line extension preserves signed-permutation invariance. Assigning the spectral weight zero on rank deficiency also preserves the joint invariance required by Theorem~\ref{thm:weighted}.

\section{The regularized spectral intensity}\label{sec:intensity}

Choose a measurable orthonormal eigenbasis of $C^\top C$ as a function of that matrix alone, and form the corresponding left singular vectors $u_k=Cq_k/s_k(C)$. On full rank these obey $u_k(PC)=Pu_k(C)$, including with a fixed measurable choice inside repeated eigenspaces. This convention makes all coordinate-flatness and truncated spectral conditions below invariant under $C\mapsto PC$.

Put $\ell=\lfloor\sqrt{\log n}\rfloor$ and
\begin{equation}\label{eq:chi}
\chi_C(X)^2=\sum_{k=1}^{\ell}\frac{\ip{u_k}{X}^2}{s_k(C)^2}.
\end{equation}
For $G_1,\ldots,G_\ell$ iid $N(0,1)$, let
\begin{equation}\label{eq:Y}
Y_C=\left(\frac1n\sum_{k=1}^{\ell}\frac{G_k^2}{s_k(C)^2}\right)^{1/2}.
\end{equation}
We shall use the SSS deleted-matrix event $\good$ specified in Proposition~\ref{prop:sss-input}. In particular it is invariant and implies $s_1(C)\ge n^{-1/2}(\log n)^{-3}$. Define
\begin{equation}\label{eq:lambda}
\Lambda_n(C)=\1_{\good}\sqrt{\frac2\pi}\,\E_GY_C.
\end{equation}
This definition sets $\Lambda_n=0$ off $\good$, including on rank deficiency. It is independent of the rare-event threshold and of the preconditioner.

\begin{proposition}[Subcritical intensity moments]\label{prop:lp}
For every $0<p<2$,
\begin{equation}\label{eq:lp}
\sup_{n\ge n_0}\E\Lambda_n(C)^p<\infty.
\end{equation}
The same conclusion holds if $\good$ is replaced by any event contained in $\{s_1(C)\ge n^{-1/2}(\log n)^{-A}\}$, for fixed $A>0$.
\end{proposition}

\begin{proof}
It suffices to consider $1<p<2$. Jensen's inequality and subadditivity of the power $p/2$ give
\begin{equation}\label{eq:lp-sum}
\Lambda_n^p\le C_p n^{-p/2}\sum_{k=1}^{\ell}
\1_{\good}s_k(C)^{-p}.
\end{equation}
We claim
\begin{equation}\label{eq:inverse-moment}
\E[\1_{\good}s_k(C)^{-p}]\le C_p n^{p/2}k^{-p},
\qquad 1\le k\le\ell.
\end{equation}
Let $C_k$ consist of the first $n-k$ columns of $C$. Eigenvalue interlacing for the principal submatrix $C_k^\top C_k$ of $C^\top C$ gives $s_1(C_k)\le s_k(C)$. By the rectangular estimate \cite[Theorem~1.1]{RVRect},
\[
\Prob\!\left(s_1(C_k)\le u\bigl(\sqrt n-\sqrt{n-k-1}\bigr)\right)
\le (Cu)^{k+1}+Ce^{-cn}.
\]
Since $\sqrt n-\sqrt{n-k-1}\ge k/(2\sqrt n)$, the variable $Z_k=\sqrt n\,s_k(C)/k$ therefore satisfies
\begin{equation}\label{eq:rect-k-tail}
\Prob(Z_k\le u)\le(C'u)^{k+1}+Ce^{-cn}.
\end{equation}
On $\good$ it obeys $Z_k\ge a_{n,k}:=(\log n)^{-3}/k$. Fix $u_*>0$ such that $C'u_*<1/2$. For sufficiently large $n$, $a_{n,k}<u_*$ uniformly for $1\le k\le\ell$, and tail integration yields
\[
\E[\1_{\good}Z_k^{-p}]
\le u_*^{-p}+p\int_{a_{n,k}}^{u_*}
u^{-p-1}(C'u)^{k+1}\,du
+Ce^{-cn}a_{n,k}^{-p}\le C_p.
\]
Here $k+1\ge2>p$, and the last term is uniformly bounded for $k\le\ell$. This proves \eqref{eq:inverse-moment}. Summing \eqref{eq:lp-sum} gives $C_p\sum_{k\ge1}k^{-p}<\infty$. Smaller positive $p$ follow by monotonicity of moments. The argument with a general fixed logarithmic exponent is identical.
\end{proof}

\begin{remark}[Why second moments are unnecessary]
The raw Gaussian deficiency-one intensity has a critical second moment: the lower tail of $\sqrt n\,s_1(C)$ has quadratic order, while the intensity is at least a constant times $(\sqrt n\,s_1(C))^{-1}$. The proof uses a fixed $p\in(1,2)$ and does not require a second-moment estimate for the intensity.
\end{remark}

\begin{lemma}[Gaussian linearization]\label{lem:gauss-linear}
For full-rank $C$, define
\[
p_C(s)=\Prob_{G_0,G}(|G_0|\le sY_C),
\]
where $G_0$ is independent of $G_1,\ldots,G_\ell$. On $\good$, for $s\ge0$,
\begin{equation}\label{eq:gauss-linear}
|p_C(s)-s\Lambda_n(C)|
\le Cs^3(\log n)^9\ell^{3/2},
\qquad p_C(s)\le s\Lambda_n(C).
\end{equation}
\end{lemma}

\begin{proof}
For every $x\ge0$, $0\le\sqrt{2/\pi}\,x-\Prob(|G_0|\le x)\le Cx^3$. On $\good$, $Y_C\le(\log n)^3\norm G_2$, and $\E\norm G_2^3\le C\ell^{3/2}$. Conditioning on $Y_C$ proves the claim.
\end{proof}

\section{A marked submicroscopic transfer}\label{sec:transfer}

We first state the precise unmarked facts used from SSS. Keeping the full regularity event separate from the deleted-matrix event is essential.

\begin{proposition}[Unmarked SSS inputs]\label{prop:sss-input}
There exist positive constants $b,\beta,\kappa$ and $\delta_n=(\log n)^{-\kappa}$, together with a full event $\reg=\reg(C,X)$ and a deleted-matrix event $\good=\good(C)$, with the following properties. The latter event is signed-permutation invariant, contains only full-rank matrices, and implies $s_1(C)\ge n^{-1/2}(\log n)^{-3}$. Uniformly for $e^{-bn}\le s\le n^{-\beta}$,
\begin{equation}\label{eq:sss-bad}
\Prob\bigl(S_s\cap(\reg^c\cup\good^c)\bigr)
\le C\delta_ns+Ce^{-bn}.
\end{equation}
On $\reg$,
\begin{equation}\label{eq:sss-trunc}
\chi_C(X)\le\Xi_C(X)\le(1+C\delta_n)\chi_C(X).
\end{equation}
For every fixed $C\in\good$ and every such $s$,
\begin{equation}\label{eq:sss-replace}
\begin{split}
p_C((1-\delta_n)s)-C\delta_ns
&\le\Prob_X(|v^\top X|\le sn^{-1/2}\chi_C(X))\\
&\le p_C((1+\delta_n)s)+C\delta_ns.
\end{split}
\end{equation}
Finally, for the iid square matrix,
\begin{equation}\label{eq:sss-scalar}
\Prob(S_s)=(1+o(1))s+O(e^{-bn})
\end{equation}
uniformly on any sufficiently small subwindow of the displayed range.
\end{proposition}

\begin{proof}[Source and conventions]
Apply \cite{SSS} to $A^\top$. Their Lemmas~4.2 and 8.2 give \eqref{eq:sss-bad}. The deterministic comparison \eqref{eq:sss-trunc} follows from equation~(27) in the proof of their Lemma~4.3. Their Lemma~8.1, equations~(49)--(50), gives \eqref{eq:sss-replace}, and their Corollary~1.3 gives \eqref{eq:sss-scalar}. The replacement estimate holds for each fixed deleted matrix in their regularity event. Decrease $b$ and $\kappa$ and increase $\beta$ if needed, choosing the displayed closed threshold interval strictly inside the source's open intervals. Take $b$ smaller also than the exponential constants in \eqref{eq:unmarked-linear} and \eqref{eq:residual}.

The event $\good$ is their $\mathcal E^*$ from Section~8.1. The singular-vector convention preceding \eqref{eq:chi} makes its spectral and coordinate-flatness conditions invariant. Its LCD condition is invariant as well, since signed permutations preserve the integer lattice and Euclidean distances. The event $\reg$ includes the fresh-vector projection conditions from their Section~4. These conditions depend on $X$, so $\reg$ must be kept distinct from the deleted-matrix event $\good$.
\end{proof}

\begin{theorem}[Transfer for orbit-bounded marks]\label{thm:transfer}
Let $a=a(C)\ge0$ be measurable, set to zero on rank deficiency, and suppose
\begin{equation}\label{eq:orbit-assumption}
\E_{P\in\Gn}a(PC)^2\le K_a
\end{equation}
for every full-rank $C$, where $K_a<\infty$ is fixed. Then there exist $0<c_0<c_1$, $D>1$, and a deterministic sequence $\eta_n\to0$, depending only on the entry law and $K_a$, such that
\begin{equation}\label{eq:transfer}
\left|\Prob(T_{ta})-t\E[a(C)\Lambda_n(C)]\right|
\le Ct\eta_n+Ce^{-c_1n},
\qquad e^{-c_0n}\le t\le n^{-D}.
\end{equation}
The estimate is uniform over all marks satisfying \eqref{eq:orbit-assumption}. No concentration of $a$ around one is required.
\end{theorem}

\begin{proof}
We may suppose $K_a\ge1$. Write $\delta=\delta_n$ from Proposition~\ref{prop:sss-input}, choose
\begin{equation}\label{eq:L-choice}
L=L_n=\min\{\delta^{-1/2},\log n\},
\qquad h=e^{-b'n},
\end{equation}
where $0<b'<b$ is fixed, and put $J=\{h/t\le a\le L\}$. Choose $c_0<c_1<b'$ and $D>\max\{\beta+1,2\}$. These choices leave fixed slack for all threshold factors $1+O(\delta)$ below. All assertions are for sufficiently large $n$.

\smallskip\noindent\emph{Removing large and very small marks.}
Lemma~\ref{lem:orbit-tail}, applied to $T_s$ and \eqref{eq:unmarked-linear}, yields
\begin{equation}\label{eq:transfer-tail}
\Prob(T_{ta}\cap\{a>L\})\le C t/L+Ce^{-bn}.
\end{equation}
By monotonicity and Lemma~\ref{lem:inclusion},
\begin{equation}\label{eq:transfer-low}
\Prob(T_{ta}\cap\{ta<h\})\le\Prob(T_h)\le Ch+Ce^{-bn}.
\end{equation}
For the coefficient, invariance of $\Lambda_n$, Lemma~\ref{lem:orbit-tail}, and Proposition~\ref{prop:lp} give
\begin{equation}\label{eq:coefficient-cutoffs}
t\E[a\Lambda_n\1_{J^c}]\le Ct/L+Ch.
\end{equation}
Indeed the large-mark part is bounded by $tK_a\E\Lambda_n/L$, and the small-mark part by $h\E\Lambda_n$.

\smallskip\noindent\emph{Unmarked regularity at the enlarged threshold.}
On $\{a\le L\}$, $T_{ta}\subset T_{tL}\subset S_{tL}$. Consequently,
\begin{equation}\label{eq:enlarged-regularity}
\Prob(T_{ta}\cap J\cap(\reg^c\cup\good^c))
\le C\delta tL+Ce^{-bn}.
\end{equation}
Only the unmarked estimate \eqref{eq:sss-bad} is used. Since $tL\le n^{-\beta}$ and $tL\ge e^{-bn}$, its range is respected.

Let
\[
q_X(s;J)=\E_C\1_{\good\cap J}
\Prob_X(|v^\top X|\le s a(C)n^{-1/2}\chi_C(X)).
\]
The inequality $\chi_C\le\Xi_C$ gives $q_X(t;J)\le\Prob(T_{ta})$. Conversely, \eqref{eq:sss-trunc}, \eqref{eq:transfer-tail}--\eqref{eq:enlarged-regularity} give, for a fixed constant $C_\star>0$,
\begin{equation}\label{eq:q-sandwich}
\begin{split}
q_X(t;J)\le\Prob(T_{ta})
\le q_X((1+C_\star\delta)t;J)
+Ct(\delta L+L^{-1})+Ch+Ce^{-bn}.
\end{split}
\end{equation}

\smallskip\noindent\emph{Conditional replacement and linearization.}
For $C\in\good\cap J$, the effective threshold $ta(C)$ lies between $h$ and $tL$. The fixed endpoint slack also permits every threshold $\alpha ta(C)$ with $1\le\alpha\le1+C_\star\delta$, so \eqref{eq:sss-replace} applies pointwise throughout this interval. Its averaged error is bounded by $C\delta t\E a\le C\delta t\sqrt{K_a}$. Lemma~\ref{lem:gauss-linear}, used at thresholds $(1\pm\delta)\alpha ta(C)$, and $a^3\1_J\le L^2a$ show that
\begin{equation}\label{eq:q-linear}
\sup_{1\le\alpha\le1+C_\star\delta}
\left|q_X(\alpha t;J)-t\E[a\Lambda_n\1_J]\right|
\le C\delta t+Ct^3L^2(\log n)^9\ell^{3/2}.
\end{equation}
The threshold displacement costs $O(\delta t)$ because
\[
\E[a\Lambda_n]=\E[\Lambda_n\E_Pa(PC)]
\le\sqrt{K_a}\E\Lambda_n\le C.
\]
Combining \eqref{eq:coefficient-cutoffs}, \eqref{eq:q-sandwich}, and \eqref{eq:q-linear} proves \eqref{eq:transfer}, for example with
\[
\eta_n=C\left(\delta L+L^{-1}+
n^{-2D}L^2(\log n)^9\ell^{3/2}\right)\to0.
\]
The terms $h$ and $e^{-bn}$ are $O(e^{-c_1n})$. None of the errors requires $\eta_n=o(t)$.
\end{proof}

\begin{remark}[Why the proof preserves exponential precision]
The truncation uses an orbitwise bound inside $T_{ta}$, rather than the unconditional probability of $\{a>L\}$. Thus its cost is $t/L$. The bounded part loses $t\delta_nL$. Taking $L\to\infty$ slowly makes both errors $o(t)$ even for $t=e^{-\alpha n}$. There is no additional requirement on the rank or condition number of the preconditioner.
\end{remark}

\begin{proposition}[Normalization of the intensity]\label{prop:normalization}
The intensity in \eqref{eq:lambda} satisfies $\E\Lambda_n(C)=1+o(1)$.
\end{proposition}

\begin{proof}
Choose a fixed $D'>D$ sufficiently large that $t=n^{-D'}$ and its bounded multiplicative perturbations lie inside both the transfer window and the scalar subwindow in \eqref{eq:sss-scalar}. Apply Theorem~\ref{thm:transfer} to the mark equal to one on full rank to get
\[
\Prob(T_t)=t\E\Lambda_n+o(t).
\]
It remains to compare $T_t$ with $S_t$. For $B=I/\sqrt n$, the rank-one norm in \eqref{eq:inverse-sandwich} is at least $1/t$ precisely on $T_t$, and $D_C\le1/(\sqrt n\,s_1(C))$. If $\smin(A)\le t/((1+\gamma)\sqrt n)$ and $D_C\le\gamma/t$, then \eqref{eq:inverse-sandwich} implies $T_t$, away from square singularity. Hence
\[
\Prob(S_{t/(1+\gamma)})-C t^2/\gamma^2-Ce^{-cn}
\le\Prob(T_t)\le\Prob(S_t).
\]
Set $\gamma=t^{1/3}$ and use \eqref{eq:sss-scalar}. Both bounds equal $t(1+o(1))$. This proves the normalization without using the anisotropic theorem.
\end{proof}

\section{Proof of the main theorem}\label{sec:main-proof}

We first prove the inverse form, which also yields the rectangular-input applications.

\begin{theorem}[Universal inverse amplification]\label{thm:inverse-main}
Under \eqref{eq:entry}, for every deterministic $B\in\R^{n\times m}$ with $\norm B_{\HS}=1$, the failure probability
\[
p_n(\varepsilon,B)=\Prob(\norm{A_n^{-1}B}_{\op}\ge\varepsilon^{-1})
\]
satisfies
\begin{equation}\label{eq:inverse-main}
|p_n(\varepsilon,B)-\varepsilon\Theta_n(B)|
\le C\varepsilon\left[r_n+
(\varepsilon n\norm B_{\op}^2)^{1/3}\right]+Ce^{-c_1n}
\end{equation}
in a window \eqref{eq:window}, uniformly in $m$ and $B$. The inverse norm is assigned $+\infty$ on square singularity.
\end{theorem}

\begin{proof}
The mark $a_B$ satisfies \eqref{eq:orbit-assumption} with $K_a=1$. By Theorem~\ref{thm:transfer},
\[
\Prob(T_{ta_B})=t\E[a_B\Lambda_n]+o(t)+O(e^{-cn})
\]
uniformly in $B$ in the transfer window. The normal line satisfies the coordinate input of Theorem~\ref{thm:weighted} by Lemma~\ref{lem:coordinate}. The invariant intensity has uniformly bounded $L^{4/3}$ norm by Proposition~\ref{prop:lp}. Therefore Theorem~\ref{thm:weighted} and Proposition~\ref{prop:normalization} give
\begin{equation}\label{eq:coefficient-identification}
\E[a_B\Lambda_n]=\Theta_n(B)+o(1)
\end{equation}
uniformly over $m,B$. Its boundedness also follows directly from orbit averaging.

Apply \eqref{eq:marked-sandwich} at thresholds $\varepsilon/(1\pm\delta)$. Choose the main window strictly inside the transfer window. The threshold displacements cost $O(\varepsilon\delta)$, while \eqref{eq:residual} costs $O(\varepsilon^2n\norm B_{\op}^2/\delta^2)+O(e^{-cn})$. Setting
\[
\delta=(\varepsilon n\norm B_{\op}^2)^{1/3}
\]
balances these errors. It is less than $1/4$ uniformly for large $n$ once $C_0>1$. Gather the uniform vanishing coefficient and transfer errors into $r_n$. Square and deleted singularity probabilities are absorbed into the exponential remainder, decreasing $c_1$ and then $c_0$ if necessary so that $c_0<c_1$. This proves \eqref{eq:inverse-main}.
\end{proof}

\begin{proof}[Proof of Theorem~\ref{thm:main}]
For nonsingular $A_n$, with $B$ as in \eqref{eq:normalization},
\[
\smin(MA_n)\le\varepsilon/\norm{M^{-1}}_{\HS}
\quad\Longleftrightarrow\quad
\norm{A_n^{-1}B}_{\op}\ge\varepsilon^{-1}.
\]
On singularity both events count failure. Theorem~\ref{thm:inverse-main} proves \eqref{eq:main-quant}. Uniform negligibility of its errors and the coefficient lower bound proved next give \eqref{eq:main-relative}.
\end{proof}

\section{The Gaussian coefficient and anisotropy profiles}\label{sec:profiles}

\begin{proposition}[Bounds and diffuse bias]\label{prop:theta}
Let $H\succeq0$ with $\Tr H=1$. For the functional $\vartheta_n$ defined in \eqref{eq:theta},
\begin{equation}\label{eq:theta-bounds}
\sqrt{\frac2\pi}\le\sqrt n\,\E|U_1|\le\vartheta_n(H)\le1,
\end{equation}
and
\begin{equation}\label{eq:theta-bias}
0\le1-\vartheta_n(H)
\le\frac{n}{n+2}\left(\Tr H^2-\frac1n\right).
\end{equation}
For a sequence $H_n$, $\vartheta_n(H_n)\to1$ if and only if $\Tr H_n^2\to0$.
\end{proposition}

\begin{proof}
The upper bound is Jensen's inequality. If $\lambda_j$ are the eigenvalues of $H$, concavity of the square root gives pointwise
\[
\sqrt{\sum_j\lambda_j nU_j^2}\ge\sum_j\lambda_j\sqrt n|U_j|.
\]
Taking expectations gives the middle lower bound. If $g=RU$ is the Gaussian polar decomposition, $R$ and $U$ are independent and $\E R\le\sqrt n$. Thus $\sqrt n\E|U_1|=\sqrt n\E|g_1|/\E R\ge\sqrt{2/\pi}$.

Set $S=nU^\top HU$. The sphere moment formulas give
\[
\E S=1,\qquad
\operatorname{Var}(S)=\frac{2n}{n+2}\left(\Tr H^2-\frac1n\right).
\]
The exact identity $1-\E\sqrt S=\tfrac12\E(\sqrt S-1)^2$, together with $(\sqrt S-1)^2\le(S-1)^2$, proves \eqref{eq:theta-bias}.

For necessity, $\vartheta_n(H_n)\to1$ implies $\sqrt S\to1$ in $L^2$, hence $S\to1$ in probability. By convexity,
\[
\E S^4\le\sum_j\lambda_j\E(nU_j^2)^4
=\frac{105n^4}{n(n+2)(n+4)(n+6)}\le105.
\]
Thus $S^2$ is uniformly integrable and $\operatorname{Var}(S)\to0$. The variance formula implies $\Tr H_n^2\to0$. Sufficiency follows from \eqref{eq:theta-bias}.
\end{proof}

\begin{proof}[Proof of Corollary~\ref{cor:diffuse}]
Combine \eqref{eq:theta-bias} with Theorem~\ref{thm:main}.
\end{proof}

\begin{remark}[Effective ranks]
For $\Tr H=1$ and $H=BB^\top$,
\[
\norm B_{\op}^4\le\Tr H^2\le\norm B_{\op}^2.
\]
Thus divergence of $1/\Tr H^2$ and divergence of the operator stable rank $1/\norm B_{\op}^2$ are equivalent qualitative conditions. They encode different quantitative information. The coefficient bias in \eqref{eq:theta-bias} is quadratic in $\Delta_4$, whereas the inverse residual is governed by $\norm B_{\op}$.
\end{remark}

\begin{proposition}[General spectral profiles]\label{prop:profiles}
Let the eigenvalues of $H_n$ be $\lambda_{n,1}\ge\cdots\ge\lambda_{n,n}\ge0$, with sum one. Suppose $\lambda_{n,j}\to\lambda_j$ for every fixed $j$, and put $d=1-\sum_{j\ge1}\lambda_j\ge0$. Then
\begin{equation}\label{eq:profile}
\vartheta_n(H_n)\longrightarrow
\E\sqrt{d+\sum_{j\ge1}\lambda_jg_j^2}.
\end{equation}
The series is almost surely finite. In particular, for $r$ equal surviving weights carrying all the mass, the limit is
\begin{equation}\label{eq:rankr}
\E\sqrt{\chi_r^2/r}
=\sqrt{\frac2r}\frac{\Gamma((r+1)/2)}{\Gamma(r/2)}.
\end{equation}
\end{proposition}

\begin{proof}
Gaussian polar decomposition gives the exact identity
\begin{equation}\label{eq:polar-coeff}
\vartheta_n(H_n)=\frac{\sqrt n}{\E\norm{g^{(n)}}_2}
\E\sqrt{\sum_{j=1}^n\lambda_{n,j}g_j^2}.
\end{equation}
The prefactor tends to one. For fixed $r$, the first $r$ weighted terms converge in distribution. The centered remainder has variance
\[
2\sum_{j>r}\lambda_{n,j}^2\le2\lambda_{n,r+1}.
\]
Its mean tends to $1-\sum_{j\le r}\lambda_j$. Taking $n\to\infty$ first and then $r\to\infty$ identifies the limit of the quadratic sum as $d+\sum_{j\ge1}\lambda_jg_j^2$, because $\lambda_{r+1}\to0$. The series has finite expectation. The square roots have second moment one and hence are uniformly integrable, proving \eqref{eq:profile}. Formula~\eqref{eq:rankr} follows by integrating the chi density.
\end{proof}

\begin{remark}[Finite-dimensional normalization]
The coefficient $\Theta_n(B)=\vartheta_n(BB^\top)$ describes the dimension-asymptotic lower tail in Theorem~\ref{thm:main}. Exact fixed-dimensional coefficients can contain additional factors tending to one, as illustrated by the tall Gaussian model in Proposition~\ref{prop:tall-gaussian-coefficient}.
\end{remark}

\section{Applications}\label{sec:applications}

\subsection{Deterministic inputs and ellipsoidal amplification}
\begin{corollary}\label{cor:inputs}
For any deterministic nonzero $K_n\in\R^{n\times m_n}$, uniformly in \eqref{eq:window},
\begin{equation}\label{eq:inputs}
\Prob\!\left(\norm{A_n^{-1}K_n}_{\op}
\ge\frac{\norm{K_n}_{\HS}}{\varepsilon}\right)
=\varepsilon\left[\Theta_n\!\left(\frac{K_n}{\norm{K_n}_{\HS}}\right)+o(1)\right]
+O(e^{-c_1n}).
\end{equation}
For one nonzero deterministic right-hand side $b_n$,
\begin{equation}\label{eq:rhs}
\Prob\!\left(\norm{A_n^{-1}b_n}_2\ge
\frac{\norm{b_n}_2}{\varepsilon}\right)
=\left(\sqrt{\frac2\pi}+o(1)\right)\varepsilon+O(e^{-c_1n}).
\end{equation}
\end{corollary}

\begin{proof}
Normalize $K_n$ and apply Theorem~\ref{thm:inverse-main}. For a single column, use Proposition~\ref{prop:profiles} with one weight equal to one.
\end{proof}

The norm in \eqref{eq:inputs} is the largest amplification of an input in the ellipsoid $K_n\mathbb B_2^{m_n}$. Taking $K_n=I$ gives worst-case amplification, whereas choosing $K_n$ with orthonormal columns gives amplification restricted to a prescribed input subspace. The theorem also covers rank-deficient input maps.

\subsection{Square sample covariance with arbitrary population anisotropy}
\begin{corollary}\label{cor:covariance}
Let $\Sigma_n$ be deterministic positive definite and set
\[
S_n=\frac1n\Sigma_n^{1/2}A_nA_n^\top\Sigma_n^{1/2},
\qquad B_n=\frac{\Sigma_n^{-1/2}}{\sqrt{\Tr\Sigma_n^{-1}}}.
\]
Then uniformly in \eqref{eq:window},
\begin{equation}\label{eq:covariance}
\Prob\!\left(\lambda_{\min}(S_n)\le
\frac{\varepsilon^2}{n\Tr\Sigma_n^{-1}}\right)
=\varepsilon[\Theta_n(B_n)+o(1)]+O(e^{-c_1n}).
\end{equation}
\end{corollary}

\begin{proof}
Use $\lambda_{\min}(S_n)=n^{-1}\smin(\Sigma_n^{1/2}A_n)^2$ and apply Theorem~\ref{thm:main}.
\end{proof}

The conclusion holds for arbitrary positive-definite $\Sigma_n$, with symmetric iid latent entries and thresholds in the stated shrinking window.

\subsection{Independent random environments}
\begin{corollary}[Conditional universality]\label{cor:environment}
Let $\mathcal E_n$ be a sigma-field independent of $A_n$, and let $M_n$ be an almost surely invertible $\mathcal E_n$-measurable matrix. Put $B_n=M_n^{-1}/\norm{M_n^{-1}}_{\HS}$. For deterministic $\varepsilon$ in \eqref{eq:window}, almost surely,
\begin{equation}\label{eq:quenched}
\Prob\!\left(\smin(M_nA_n)\le
\frac{\varepsilon}{\norm{M_n^{-1}}_{\HS}}\,\middle|\,\mathcal E_n\right)
=\varepsilon[\Theta_n(B_n)+O(\widetilde r_n)]+O(e^{-c_1n}),
\end{equation}
where $\widetilde r_n\to0$ and all error constants are deterministic and uniform in the environment.
\end{corollary}

\begin{proof}
Condition on $\mathcal E_n$ and use the uniform bound of Theorem~\ref{thm:main}. The residual is bounded by $(\varepsilon n)^{1/3}$ and can be included in $\widetilde r_n$.
\end{proof}

For two independent fresh iid matrices, also independent of the same environment, let $E_{n,1},E_{n,2}$ be the corresponding failures with the same $M_n$ and threshold. Conditional independence and \eqref{eq:quenched} imply
\begin{equation}\label{eq:environment-cov}
\varepsilon^{-2}\operatorname{Cov}(\1_{E_{n,1}},\1_{E_{n,2}})
=\operatorname{Var}(\Theta_n(B_n))+o(1).
\end{equation}
Indeed the covariance equals the variance of the conditional failure probability, and its ratio to $\varepsilon$ differs uniformly from $\Theta_n(B_n)$ by $o(1)$. Thus concentration of the coefficient is sufficient for asymptotic decorrelation in this conditional-independence setting.

\subsection{Independent exponential families}
\begin{corollary}[Independent extremes]\label{cor:extremes}
Let $N_n=e^{hn+o(n)}$ with $0<h<\min\{c_0,c_1\}$. Let $A_{n,j}$, $1\le j\le N_n$, be independent copies of $A_n$, and let $M_{n,j}$ be deterministic invertible matrices. Define
\[
Z_{n,j}=\norm{M_{n,j}^{-1}}_{\HS}\smin(M_{n,j}A_{n,j}),\qquad
\theta_{n,j}=\Theta_n\!\left(\frac{M_{n,j}^{-1}}{\norm{M_{n,j}^{-1}}_{\HS}}\right).
\]
If $N_n^{-1}\sum_{j\le N_n}\theta_{n,j}\to\bar\theta$, then for each $z>0$,
\begin{equation}\label{eq:poisson}
\sum_{j=1}^{N_n}\1_{\{Z_{n,j}\le z/N_n\}}
\ \Longrightarrow\ \operatorname{Poisson}(z\bar\theta),
\end{equation}
and $N_n\min_{j\le N_n}Z_{n,j}$ converges in distribution to an exponential random variable of rate $\bar\theta$.
\end{corollary}

\begin{proof}
The threshold $z/N_n$ lies in \eqref{eq:window} for large $n$. Theorem~\ref{thm:main}, uniformly in $j$, gives probabilities $p_{n,j}$ with $\sum_jp_{n,j}\to z\bar\theta$ and $\max_jp_{n,j}\to0$. The accumulated exponential floor tends to zero because $N_ne^{-c_1n}\to0$. Independence and the product of the Bernoulli probability-generating functions prove \eqref{eq:poisson}. The zero-count probability gives the minimum law.
\end{proof}

\section{Scope of the transfer principle}\label{sec:scope}

Theorems~\ref{thm:weighted} and \ref{thm:transfer} have distinct roles. The former permits any invariant $L^p$ weight and requires no quantitative convergence rate relative to the rare-event threshold. The latter carries a measurable mark through an exponentially accurate geometric comparison using an orbitwise second-moment bound. Their combination applies whenever a rare-event problem produces these two structures.

The symmetry assumption is used for signed-permutation invariance. Removing it would require a replacement for the orbit identities, not merely a change in the scalar small-ball theorem. Dependence between a preconditioner and the deleted matrix also needs separate analysis. For example, with $v=v(C)$ and $\delta_n\downarrow0$, the random positive-definite matrix
\[
H(C)=\frac{I-vv^\top+\delta_nvv^\top}{n-1+\delta_n}
\]
has trace one and diffuse squared spectral mass, but $nv^\top H(C)v=n\delta_n/(n-1+\delta_n)\to0$. Thus diffuseness alone does not control an adapted mark.

For overlapping minors, an application would require both the conditional law of the fresh matrix and estimates for joint failures. General orthogonal rotations of non-Gaussian iid entries do not preserve independence, so these points must be checked separately. In particular, extending Corollary~\ref{cor:extremes} to dependent families requires joint estimates beyond the independent-environment results proved here.

\section{Conjectures and further directions}\label{sec:directions}

Theorem~\ref{thm:main} expresses the leading lower-tail coefficient as a Haar functional of the normalized inverse ellipsoid. This geometric description suggests extensions both to larger hard-edge thresholds and to rectangular matrices. At larger thresholds, the natural comparison object is the full Gaussian least-singular-value distribution. For tall rectangular matrices, volumes of higher-dimensional orthogonal projections of the inverse ellipsoid suggest analogues of the square mean-width coefficient. The conjectures below formulate these extensions when the row and column dimensions differ by a fixed amount.

Correlated and spiked Gaussian hard-edge laws \cite{HHN,WG,RR} and quantitative universality results for unpreconditioned iid matrices \cite{TV,Wang} provide natural antecedents. Here the proposed comparison with the Gaussian model is uniform over all deterministic invertible preconditioners. For the rectangular coefficient asymptotic, this uniformity also concerns the Gaussian small-threshold remainder when the inverse profile approaches rank deficiency.

\subsection{Universality from rare tails to the ordinary hard edge}

Fix an integer $q\ge1$. Let $A^+_{n,q}$ be an $n\times(n+q-1)$ matrix with iid entries of law $\xi$, and let $G^+_{n,q}$ have iid standard real Gaussian entries. Here $\smin$ denotes the least of the $n$ singular values, and equals zero on rank deficiency. For $M\in GL_n(\R)$, retain the normalization $h(M)$ from \eqref{eq:normalization} and write
\[
F^\xi_{n,q,M}(t)=\Prob(h(M)\smin(MA^+_{n,q})\le t),
\]
and define $F^G_{n,q,M}$ analogously.

\begin{conjecture}[Uniform Gaussian comparison across hard-edge scales]\label{conj:full-hardedge}
Let $\xi$ be a fixed real, centered, variance-one, subgaussian random variable. Symmetry is not assumed. For each fixed $q\ge1$, there exists $c_{\xi,q}>0$ such that, for every fixed $T>0$,
\begin{equation}\label{eq:conj-full}
\sup_{M\in GL_n(\R)}
\ \sup_{e^{-c_{\xi,q}n}\le t\le T}
\left|\frac{F^\xi_{n,q,M}(t)}{F^G_{n,q,M}(t)}-1\right|
\longrightarrow0.
\end{equation}
\end{conjecture}

For $M=I$ and $q=1$, the comparison in Conjecture~\ref{conj:full-hardedge} already follows from \cite[Theorem~1.2]{SSS}, including for nonsymmetric entry laws, after choosing the exponential endpoint sufficiently small. For arbitrary $M$, Theorem~\ref{thm:main} establishes the symmetric square case only in the window \eqref{eq:window}. The conjecture extends this anisotropic comparison to ordinary hard-edge thresholds, nonsymmetric entries, and every fixed excess of columns. An additive $o(1)$ comparison at ordinary thresholds would still leave an intermediate shrinking-threshold range unless its error bound overlaps with that of the present theorem.

The Gaussian reference must retain the full inverse profile. In the square case $M=I$, its limiting distribution is $1-e^{-t-t^2/2}$ \cite{Edelman,TV}. For a rank-one input, the reciprocal Gaussian amplification is exactly $|g|$. More generally, the finite-rank identity in Proposition~\ref{prop:gaussian-benchmarks} below gives explicit benchmarks, while spiked Gaussian hard-edge laws provide further antecedents \cite{RR}. The role of $\Theta_n$ in Theorem~\ref{thm:main} is to describe the first-order rare-tail behavior.

The exponential endpoint must depend on the entry law and on $q$. For example, for $\xi=(\operatorname{Bernoulli}(p)-p)/\sqrt{p(1-p)}$, coincidence of two rows of $A^+_{n,q}$ has probability $[p^2+(1-p)^2]^{n+q-1}$ and forces rank deficiency. Such atoms must be negligible compared with the Gaussian lower tail. In contrast with the tall model below, arbitrary inverse anisotropy does not make the wide Gaussian coefficient vanish.

For this wide model, the candidate dimension-asymptotic small-threshold coefficient is
\begin{equation}\label{eq:wide-coefficient}
\frac1{q!}\Psi_{n,q}(H),\qquad
\Psi_{n,q}(H)=\E_{U\sim\Haar(S^{n-1})}(nU^\top HU)^{q/2},
\qquad H=\frac{M^{-1}M^{-\top}}{h(M)^2}.
\end{equation}
The proposed asymptotic is a lower tail of order $t^q$ with coefficient $\Psi_{n,q}(H)/q!$ as $n\to\infty$. Exact fixed-dimensional Gaussian coefficients may include further normalization factors. At $q=1$ it is the coefficient of this paper. At $q=2$, $\Psi_{n,2}=1$ exactly. For $q>2$, the functional is convex and symmetric in the eigenvalues of $H$, so concentrating inverse singular mass increases it in the majorization order, reversing the concave ordering at $q=1$. A rank-one profile gives $\E|g|^q/q!=1/[2^{q/2}\Gamma(q/2+1)]$ in the large-$n$ limit, matching the lower-tail coefficient of $\chi_q$.

\subsection{A rectangular intrinsic-volume hierarchy}

Now let $A^-_{n,q}$ have size $n\times(n-q+1)$ and iid entries of law $\xi$, where $q\ge1$ is fixed and $n\ge q$. The multiplier still acts on the $n$ rows. For trace-one $H\succeq0$, define
\begin{equation}\label{eq:grassmann-coefficient}
\Phi_{n,q}(H)=
\E_Q\sqrt{\det(nQ^\top HQ)},
\qquad Q^\top Q=I_q,
\end{equation}
where the columns of $Q$ form a Haar-uniform orthonormal $q$-frame in $\R^n$. Thus $\Phi_{n,1}(H)=\vartheta_n(H)$.

\begin{conjecture}[Intrinsic-volume coefficients at a rectangular hard edge]\label{conj:intrinsic}
Assume \eqref{eq:entry}. For each fixed $q\ge1$, there are constants $c_0,c_1,C_0,C>0$, with $qc_0<c_1$, and a deterministic $r_{n,q}\to0$, depending only on $q$ and the entry law, such that, uniformly over $M\in GL_n(\R)$ and
$e^{-c_0n}\le\varepsilon\le n^{-C_0}$,
\begin{equation}\label{eq:conj-intrinsic}
\left|\Prob\!\left(\smin(MA^-_{n,q})
\le\frac{\varepsilon}{h(M)}\right)
-\frac{\varepsilon^q}{q!}\Phi_{n,q}(H)\right|
\le r_{n,q}\varepsilon^q+Ce^{-c_1n},
\end{equation}
where $H=M^{-1}M^{-\top}/h(M)^2$ and $\smin$ is zero on rank deficiency.
\end{conjecture}

The coefficients have an intrinsic-volume interpretation. Let $\kappa_j$ be the volume of the unit ball in $\R^j$, and let $V_q$ denote the intrinsic volume with the Steiner normalization. Kubota's formula gives
\begin{equation}\label{eq:intrinsic-volume}
\Phi_{n,q}(H)=
\frac{n^{q/2}\kappa_{n-q}}{\binom nq\kappa_n}
V_q(H^{1/2}\mathbb B_2^n).
\end{equation}
See \cite{GSZ} for this normalization and ellipsoid formulas. The square coefficient is therefore the first member, corresponding to mean width. The conjecture would identify all fixed-order intrinsic volumes of the inverse ellipsoid as coefficients of real rectangular hard-edge tails.

The distinction between the two rectangular models is essential. One has $0\le\Phi_{n,q}\le1$, but for $q\ge2$ this coefficient can tend to zero. Indeed, Cauchy--Binet and Haar symmetry give
\[
\E_Q\det(nQ^\top HQ)=
\frac{n^q e_q(H)}{\binom nq}\le1,
\]
where $e_q$ is the $q$th elementary symmetric polynomial. Jensen's inequality then gives the upper bound on $\Phi_{n,q}$. If all inverse mass concentrates on $r<q$ limiting directions, the projected determinant vanishes. Under a finite-spike profile with weights $\lambda_1,\ldots,\lambda_r$ and diffuse residual mass $d$, the corresponding geometric limit is
\begin{equation}\label{eq:grassmann-profile}
\Phi_{n,q}(H_n)\longrightarrow
\E\sqrt{\det\!\left(dI_q+\sum_{j=1}^r\lambda_j g_jg_j^\top\right)},
\qquad g_j\sim N(0,I_q)\ \text{independently}.
\end{equation}
This follows by the Gaussian limit of a fixed number of Haar rows and concentration of the diffuse remainder. The preceding determinant identity supplies uniform integrability of the square roots.

Consequently \eqref{eq:conj-intrinsic} gives a relative asymptotic when $\Phi_{n,q}$ is bounded below. A relative comparison uniform over all tall preconditioners, without an additive floor, would be false for discrete entries. For example, duplicate columns give a positive rank-deficiency probability for Rademacher matrices. With $M_a=\diag(a,1,\ldots,1)$ and $q\ge2$, the Gaussian lower $(n-1)\times(n-q+1)$ block has full column rank almost surely, while $\varepsilon/h(M_a)\to0$ as $a\downarrow0$. The Gaussian failure probability therefore tends to zero, whereas the discrete singularity atom is unchanged.

\subsection{Exact Gaussian identities}

The following propositions give exact Gaussian identities that fix the normalizations in the two rectangular models.

\begin{proposition}[Gaussian benchmarks]\label{prop:gaussian-benchmarks}
Let $G^+_{n,q}$ be as above. If $V\in\R^{n\times r}$ has orthonormal columns and $D$ is a positive diagonal $r\times r$ matrix, then
\begin{equation}\label{eq:gaussian-rank-reduction}
\norm{(G^+_{n,q})^\dagger VD}_{\op}
\overset d=
\norm{(G^+_{r,q})^\dagger D}_{\op}.
\end{equation}
Here $\dagger$ denotes the Moore--Penrose inverse. In particular, the reciprocal for $r=1,D=1$ has law $\chi_q$. Moreover, for $q=2$, every $n$, every invertible $M$, and every $t\ge0$,
\begin{equation}\label{eq:gaussian-rayleigh}
F^G_{n,2,M}(t)=1-e^{-t^2/2}.
\end{equation}
\end{proposition}

\begin{proof}
Put $W=G^+_{n,q}(G^+_{n,q})^\top$. The squared norm in \eqref{eq:gaussian-rank-reduction} equals
\[
\lambda_{\max}(DV^\top W^{-1}VD).
\]
Rotate $V$ to the first $r$ coordinate vectors and condition on the remaining $n-r$ Gaussian rows. Their orthogonal complement has dimension $r+q-1$, and the projected first $r$ rows form an iid $r\times(r+q-1)$ Gaussian matrix. The Schur complement identifies the leading $r\times r$ block of $W^{-1}$ with the inverse Gram matrix of this projection, proving \eqref{eq:gaussian-rank-reduction}.

For \eqref{eq:gaussian-rayleigh}, $W=G^+_{n,2}(G^+_{n,2})^\top$ has density proportional to $e^{-\Tr W/2}$ on the positive-definite cone. The event in the survival probability is $W>t^2H$. Translation by $t^2H$ multiplies its integral by $e^{-t^2\Tr H/2}=e^{-t^2/2}$. This is also a special case of the correlated Gaussian Wishart formulas in \cite{WG}.
\end{proof}

\begin{proposition}[The tall Gaussian coefficient]\label{prop:tall-gaussian-coefficient}
Fix $n\ge2$ and $M\in GL_n(\R)$. Let $G$ be an $n\times(n-1)$ matrix with iid standard real Gaussian entries and $z$ uniform on $S^{n-1}$. With $h=h(M)$ and $H=M^{-1}M^{-\top}/h^2$,
\begin{equation}\label{eq:tall-gaussian-coefficient}
\begin{aligned}
\lim_{\varepsilon\downarrow0}\frac{1}{\varepsilon^2}
\Prob\!\left(\smin(MG)\le\frac{\varepsilon}{h}\right)
&=\frac12\left(1-\E_z\frac{z^\top H^2z}{z^\top Hz}\right)
\\
&=\frac{n-1}{2n}\Phi_{n,2}(H).
\end{aligned}
\end{equation}
\end{proposition}

\begin{proof}
Replacing $M$ by $hM$, we may assume $h=1$. Put $p=n-1$. The joint density of the ordered squared singular values $\lambda_1<\cdots<\lambda_p$ of $G$ is proportional to
\[
e^{-\sum_{i=1}^p\lambda_i/2}\prod_{i<j}(\lambda_j-\lambda_i).
\]
The substitution $X=\lambda_1$ and $Y_i=\lambda_i-\lambda_1$, $2\le i\le p$, shows that $X$ is exponential with rate $p/2$, independently of $(Y_2,\ldots,Y_p)$. Choose the left singular frame $U$ with independent uniform column signs. It is Haar and independent of these variables. Set
\[
D=\diag(0,Y_2,\ldots,Y_p),\qquad K=U^\top M^\top MU.
\]
For $p=1$, take $D=(0)$. The event $\smin(MG)^2>\varepsilon^2$ is equivalent to $XI_p+D>\varepsilon^2K^{-1}$, and hence to $X>\lambda_{\max}(\varepsilon^2K^{-1}-D)$. Thus
\[
\Prob(\smin(MG)\le\varepsilon\mid U,D)
=1-\exp\!\left[-\frac p2\lambda_{\max}(\varepsilon^2K^{-1}-D)\right].
\]
The largest eigenvalue on the right is positive, since its Rayleigh quotient at the first coordinate vector is $\varepsilon^2(K^{-1})_{11}$. The zero eigenvalue of $-D$ is simple almost surely, so
\[
\varepsilon^{-2}\lambda_{\max}(\varepsilon^2K^{-1}-D)
\longrightarrow(K^{-1})_{11}.
\]
If $z$ is a unit normal to the columns of $U$, the Schur complement yields
\[
K^{-1}=U^\top\left(H-\frac{Hzz^\top H}{z^\top Hz}\right)U.
\]
In particular, $\norm{K^{-1}}_{\op}\le\norm H_{\op}$, and the conditional probability divided by $\varepsilon^2$ is bounded by $p\norm H_{\op}/2$. Dominated convergence therefore gives the coefficient $p\E(K^{-1})_{11}/2$. Conditional on $z$, the first column of $U$ is uniform on the unit sphere of $z^\perp$. Averaging its quadratic form gives the first equality in \eqref{eq:tall-gaussian-coefficient}.

For the second equality, the ellipsoid formula in \cite[Section~2]{GSZ} gives
\[
V_2(H^{1/2}\mathbb B_2^n)
=\pi\left(\Tr H-\E_z\frac{z^\top H^2z}{z^\top Hz}\right).
\]
Now apply \eqref{eq:intrinsic-volume} with $q=2$, using $\Tr H=1$ and $\kappa_n=2\pi\kappa_{n-2}/n$.
\end{proof}

This identity determines the coefficient at fixed dimension and fixed preconditioner. A remainder bound uniform over preconditioners, as required by Conjecture~\ref{conj:intrinsic}, remains to be established.

\subsection{Further inputs required}

Conjecture~\ref{conj:full-hardedge} would require control of the bottom singular values jointly with their projections onto deterministic directions, and truncation of the remaining inverse spectrum. At ordinary thresholds several low singular modes can contribute. Removing symmetry also requires a substitute for signed-permutation averaging. For $q>1$, the rarity scale is $\varepsilon^q$, so a direct extension of the present truncation needs higher-order mark control beyond the orbitwise second-moment identity.

Conjecture~\ref{conj:intrinsic} additionally calls for comparison of a random $q$-dimensional normal space and its volume mark at matrix-dependent thresholds. Together, the conjectures concern arbitrary deterministic population anisotropy throughout the real fixed-codimension least-singular-value problem. They make no assertion about the full eigenvalue process, proportional-aspect soft edges, or growing rectangular excess. Quantitative rates for the existing theorem remain a separate question: \eqref{eq:uniform-orbit-error} pairs an $L^4$ Wasserstein bound for $\mu_{v_n}$ with the uniform $L^{4/3}$ bound for $\Lambda_n$.

\appendix
\section{The two-coordinate magnitude limit}\label{app:coordinate}

We give the fixed-coordinate consequence needed in Lemma~\ref{lem:coordinate}, using the projected-column method of Nguyen--Vu \cite[Section~6]{NV}. The sampling concentration estimate in their Theorem~6.1 is stated without dependence on the sample size. Since convergence in probability of the normalization is sufficient here, we replace that step by the finite-population variance bound below. The other inputs are their maximal-coordinate delocalization, their frame Berry--Esseen estimate, and polynomial flatness of the omitted-column complement. All references in this appendix to results of \cite{NV} use the published version.

\begin{lemma}[Sampled normalization]\label{lem:sample-normalization}
Let $n\ge4$, $2\le m\le n$, and $z_1,\ldots,z_n\in[0,L]$ with $L\ge1$ and $\sum_i z_i=n$. Let $I$ be a uniformly chosen $m$-element subset of $\{1,\ldots,n\}$ containing $1$ and $2$, and set
\[
S=\frac1m\sum_{i\in I}z_i.
\]
Then
\begin{equation}\label{eq:sample-normalization}
|\E_I S-1|\le\frac{4L}{m},
\qquad
\operatorname{Var}_I(S)\le\frac{2L}{m}.
\end{equation}
In particular, $S\to1$ in probability whenever $L/m\to0$.
\end{lemma}

\begin{proof}
Put $a_0=z_1+z_2$. The remaining $m-2$ entries are sampled without replacement from a population of size $n-2$, and
\[
\E_I(mS)-m=\frac{n-m}{n-2}(a_0-2).
\]
This gives the first bound. The finite-population variance formula, followed by $z_i^2\le Lz_i$, gives
\[
\operatorname{Var}_I(mS)
\le\frac{m-2}{n-2}\sum_{i=3}^n z_i^2
\le\frac{mnL}{n-2}\le2mL.
\]
The conclusion follows from these bounds and Chebyshev's inequality.
\end{proof}

\begin{proof}[Proof of Lemma~\ref{lem:coordinate}]
Write $\mathcal A=C_n^\top$, so $\mathcal A$ is an $(n-1)\times n$ iid matrix, and let $v$ be its unit normal on full rank. The rectangular estimate \cite[Theorem~1.1]{RVRect} makes rank deficiency exponentially unlikely. Coordinate delocalization \cite[Theorem~1.4, equation~(6)]{NV} supplies a constant $C$ such that
\begin{equation}\label{eq:coordinate-delocalization}
\Prob\!\left(\max_{1\le i\le n}n|v_i|^2>C\log n\right)=o(1).
\end{equation}

\smallskip\noindent\emph{Sampling and normalization.}
Independently of $\mathcal A$, choose a uniform $m$-element set $I$ containing $1$ and $2$, ordered with these indices first. The parameter $m$ will be a sufficiently small positive power of $n$. On the event in \eqref{eq:coordinate-delocalization} with the reverse inequality, apply Lemma~\ref{lem:sample-normalization} conditionally to $z_i=n|v_i|^2$ and $L=C\log n$. Since $L/m\to0$, it follows that
\begin{equation}\label{eq:sample-mass}
S_n:=\frac n m\sum_{i\in I}|v_i|^2\longrightarrow1
\quad\text{in probability}.
\end{equation}
This estimate will be used after the comparison below. We do not condition the retained entries on the event in \eqref{eq:coordinate-delocalization} or on a restriction on $S_n$.

\smallskip\noindent\emph{Projection and Gaussian comparison.}
Let $H$ be the span of the columns of $\mathcal A$ indexed by $I^c$. With exponentially high probability, those columns have full column rank, so $\dim H^\perp=m-1$. The flatness estimate \cite[Lemma~6.3]{NV} gives a constant $c>0$ such that, for sufficiently small polynomial choices of $m$, with probability $1-o(1)$ over $I$ and the omitted columns,
\begin{equation}\label{eq:projection-flatness}
\sup_{u\in H^\perp,\ \norm u_2=1}\norm u_\infty\le n^{-c}.
\end{equation}
Choose a fixed $\alpha>0$ with $\alpha\le\min\{c/100,1/100\}$ small enough to meet that estimate, and put
\[
m=\lfloor n^\alpha\rfloor,\qquad \rho=n^{-c/4}.
\]
The constant $c$ in the flatness estimate can be taken independently of this sufficiently small choice of $\alpha$.

Condition only on $I$ and an omitted-column realization satisfying full column rank and \eqref{eq:projection-flatness}. Choose an orthonormal basis of $H^\perp$, and let $\pi:\R^{n-1}\to\R^{m-1}$ be the corresponding coordinate projection. Thus $\pi\pi^\top=I_{m-1}$. Set $Z=\pi\mathcal A_I$, an $(m-1)\times m$ matrix. The entries of $\mathcal A_I$ are still independent copies of the original entry variable under this conditioning.

For $1\le r\le n-1$ and $1\le j\le m$, let $w_{rj}$ be the vectorization of the matrix $(\pi e_r)e_j^\top$. Then
\[
\operatorname{vec}(Z)=\sum_{r,j}\xi_{rj}w_{rj},
\qquad
\sum_{r,j}w_{rj}w_{rj}^\top=I_{m(m-1)},
\qquad
\max_{r,j}\norm{w_{rj}}_\infty\le n^{-c}.
\]
The Gaussian replacement of this normalized frame is therefore the vectorization of an iid standard Gaussian matrix $G$ of the same dimensions. For any measurable set $\Omega$ in the matrix space, the frame Berry--Esseen estimate \cite[Lemma~2.5]{NV} gives
\begin{equation}\label{eq:frame-comparison}
\begin{split}
\Prob_G(G\in\Omega\setminus\partial_\rho\Omega)-r_n^{\mathrm{BE}}
&\le\Prob(Z\in\Omega\mid I,\mathcal A_{I^c})\\
&\le\Prob_G(G\in\Omega\cup\partial_\rho\Omega)+r_n^{\mathrm{BE}},
\end{split}
\end{equation}
where $\partial_\rho\Omega$ is the $\rho$-neighborhood of the boundary in the entrywise maximum norm, and
\[
r_n^{\mathrm{BE}}\le Cm^5\rho^{-3}n^{-c}
\le Cn^{-c/5}=o(1).
\]
The constants are uniform over the conditioned realizations.

\smallskip\noindent\emph{Control of the event boundary.}
For fixed $t_1,t_2\ge0$, let $\Omega$ consist of full-row-rank matrices whose unit kernel vector $y$, defined up to sign, satisfies $\sqrt m|y_j|\le t_j$ for $j=1,2$. The Gaussian rectangular bound gives
\begin{equation}\label{eq:gaussian-gap}
\Prob_G(\smin(G)<m^{-2})\le Cm^{-3}+Ce^{-c'm},
\end{equation}
where $\smin(G)$ denotes its smallest nonzero singular value on full row rank.

Suppose $\smin(G)\ge m^{-2}$ and $G'$ differs from $G$ by at most $\rho$ in each entry. Then $\norm{G'-G}_{\op}\le m\rho=o(m^{-2})$, so $G'$ has full row rank. If $y'$ is its unit kernel vector, projection onto the row space of $G$ gives
\[
\norm{\operatorname{Proj}_{\operatorname{row}(G)}y'}_2
\le\frac{\norm{Gy'}_2}{\smin(G)}
\le m^3\rho.
\]
After choosing the sign of a unit vector $y_G\in\ker G$, we obtain
\[
\norm{y'-y_G}_2\le Cm^3\rho,\qquad
\bigl|\sqrt m|y'_j|-\sqrt m|(y_G)_j|\bigr|
\le Cm^{7/2}\rho,\quad j=1,2.
\]
The Gaussian kernel line is Haar. Each scaled coordinate magnitude of a Haar vector in $\R^m$ has a density bounded uniformly for sufficiently large $m$. Consequently,
\[
\Prob_G(G\in\partial_\rho\Omega)
\le C m^{-3}+Ce^{-c'm}+C m^{7/2}\rho=o(1).
\]
Here $m^{7/2}\rho\le n^{-43c/200}$. Combining this estimate with \eqref{eq:frame-comparison}, and then adding the probability of omitted-column irregularity, proves
\begin{equation}\label{eq:projected-coordinates}
(\sqrt m|y_1|,\sqrt m|y_2|)
\ \Longrightarrow\ (|g_1|,|g_2|).
\end{equation}
One may define $y$ arbitrarily on failure of full row rank. When $\mathcal A$ and the omitted block both have full rank, the rank decomposition under projection implies $\rank Z=m-1$. Their exceptional probabilities are exponentially small.

\smallskip\noindent\emph{Return to the original coordinates.}
On this full-rank event, $v_I\ne0$: otherwise the full-column-rank omitted block would annihilate the nonzero vector $v_{I^c}$. Projecting $\mathcal A v=0$ shows that the kernel line of $Z$ is spanned by $v_I$, so $y=v_I/\norm{v_I}_2$ up to sign. Therefore
\[
(\sqrt n|v_1|,\sqrt n|v_2|)
=\sqrt{S_n}\,(\sqrt m|y_1|,\sqrt m|y_2|).
\]
Equations~\eqref{eq:sample-mass} and \eqref{eq:projected-coordinates}, together with Slutsky's theorem, give the required limit. Independence between $S_n$ and $y$ is unnecessary. Finally, the Haar-line extension on rank deficiency changes none of these limits.
\end{proof}

\end{document}